# Stopping strategies and gambler's ruin problem

**Theo van Uem**

School of Technology, Amsterdam University of Applied Sciences,
Weesperzijde 190, 1097 DZ Amsterdam, The Netherlands.

Email: t.j.van.uem@hva.nl

Abstract
We obtain absorption probabilities and expected time until absorption for different stopping strategies in gambler's ruin problem using the concept of multiple function barriers.

1. **Introduction**

In gambler's ruin problem a well known stopping strategy is: go on until ruin or reaching a predefined boundary. In this paper we use a more flexible concept: a player may stop in a set of predefined barriers. We use the concept of a multiple function barrier (mfb) to model our stopping strategies. In a mfb there is a positive probability s that the player will immediately stop the game and probabilities p(1-s) and q(1-s) that the player will go one step forward or backward (p+q=1). In all other states (except the absorption barrier 0) the player will go one step forward or backward with probabilities p and q.
We consider three players, A,B and C, with different stopping strategies. All players start with initial capital $i_0$ and their bet is always one unit of capital. Arriving at state 0 will always give immediate absorption in that state: ruin has happened. $M_A = \{ki_0 | k \geq 1\}$ is the set of multiple function barriers for player A. We notice that player A will stop at the very beginning of the game with probability s.
Player B has the same strategy as player A, except at time 0; the state $i_0$ is a delayed multiple function barrier for t>0. At t=0 it is a normal state (player B will always start the game, player A is a fictive player which is easier to analyze and there is a simple relationship between A and B).
Player C is more risk seeking than B: mfb's are $M_C = \{ki_0 | k \geq 2\}$; $i_0$ is a normal state, also for t>0.
In the appendix we consider the special cases s=0 and s=1.

**2. Expected number of arrivals before absorption**

We define:

$p_{ij}^{(n)}$=P(system is in state j after n steps without absorption | start in i).

The moment generating function for player X:

$X_j = X_{i_0,j}(z) = \sum_{m=0}^{\infty} p_{i_0,j}^{(m)} z^m \quad (0 < z \leq 1)$.

$\omega = \frac{p}{q}$

In the main text we suppose 0<s<1.

**Theorem 1**
The moment generating functions U for player A on the set $\{ki_0 \mid k \geq 0\}$ are:

CASE [0<z<1]: $\quad U_0 = U_0(z, \omega, s) = \omega^{-i_0} \varphi_2$

$$U_{ki_0} = U_{ki_0}(z, \omega, s) = \frac{(\tau_2^{i_0} - \tau_1^{i_0}) \varphi_2^k}{q(1-s)z(\tau_2 - \tau_1)\omega^{i_0}} \quad (k \geq 1)$$

where $\varphi_1, \varphi_2$ are solutions of $\varphi^2 - \theta\varphi + \omega^{i_0} = 0 \quad (\varphi_1 > 1 > \varphi_2 > 0)$

$$\theta = \theta(z, \omega, s) = \frac{\frac{\tau_2^{i_0} - \tau_1^{i_0}}{1-s} - 2pz(\tau_2^{i_0-1} - \tau_1^{i_0-1})}{qz(\tau_2 - \tau_1)}$$

and $\tau_1, \tau_2$ are solutions of $\quad qz\tau^2 - \tau + pz = 0 \quad (\tau_1 > 1 > \tau_2 > 0)$

CASE [z=1 and $\omega \neq 1$] $\quad \tau_1 = \max(1, \omega), \tau_2 = \min(1, \omega)$.

CASE [z=1 and $\omega = 1$] $\quad U_0 = U_0(1,1,s) = \varphi_2$

$$U_{ki_0} = U_{ki_0}(1,1,s) = \frac{2i_0 \varphi_2^k}{1-s} \quad (k \geq 1)$$

$$\varphi^2 - \theta\varphi + 1 = 0 \quad (\varphi_1 > 1 > \varphi_2 > 0)$$

$$\theta = \theta(1,1,s) = 2\left(\frac{i_0}{1-s} + 1 - i_0\right)$$

Proof

We prove the results of CASE [0<z<1] in 5 steps:

Step 1
We determine the structure of the random walk between the multiple function barriers:

$$U_{ki_0+n} = \sum_{m=1}^{\infty} p_{i_0, ki_0+n}^{(m)} z^m = \sum_{m=1}^{\infty} \left[ p_{i_0, ki_0+n-1}^{(m-1)} \cdot p + p_{i_0, ki_0+n+1}^{(m-1)} \cdot q \right] z^m = pz U_{ki_0+n-1} + qz U_{ki_0+n+1}$$

where at first glance $k \geq 0$, $n = 2, 3, \ldots, i_0 - 2$, but n=1 and n=$i_0 - 1$ also satisfy the pattern.

We get: $\quad U_{ki_0+n} = a_{k+1}\tau_1^n + b_{k+1}\tau_2^n \quad (k \geq 0, \quad n = 1, 2, \ldots, i_0 - 1)$

and $\tau_1, \tau_2$ are solutions of $\quad qz\tau^2 - \tau + pz = 0 \quad (\tau_1 > 1 > \tau_2 > 0)$

Step 2
We express the constants of step 1 in the U values of the mfb's.
On the interval $[ki_0, (k+1)i_0]$, $k \geq 1$ we have
$$U_{ki_0+1} = p(1-s)zU_{ki_0} + qzU_{ki_0+2}$$
$$U_{(k+1)i_0-1} = pzU_{(k+1)i_0-2} + q(1-s)zU_{(k+1)i_0}$$

$$a_{k+1} = \frac{(1-s)\left[\tau_2^{i_0} U_{ki_0} - U_{(k+1)i_0}\right]}{\tau_2^{i_0} - \tau_1^{i_0}}$$

$$b_{k+1} = \frac{(1-s)\left[U_{(k+1)i_0} - \tau_1^{i_0} U_{ki_0}\right]}{\tau_2^{i_0} - \tau_1^{i_0}}$$

**Step 3:**
Focussing on the mfb in $ki_0$ we get:
$$U_{ki_0} = pzU_{ki_0-1} + qzU_{ki_0+1} \quad (k > 1)$$
and after some calculation:
$$U_{(k+1)i_0} - \theta U_{ki_0} + \omega^{i_0} U_{(k-1)i_0} = 0 \quad (k > 1)$$

where
$$\theta = \theta(z, \omega, s) = \frac{\frac{\tau_2^{i_0} - \tau_1^{i_0}}{1-s} - 2pz(\tau_2^{i_0-1} - \tau_1^{i_0-1})}{qz(\tau_2 - \tau_1)}$$

and $\tau_1, \tau_2$ are solutions of $\quad qz\tau^2 - \tau + pz = 0 \quad (\tau_1 > 1 > \tau_2 > 0)$

So we have: $\quad U_{ki_0} = C_2 \varphi_2^k \quad (k \geq 1)$,
where $\varphi_1, \varphi_2$ are solutions of $\varphi^2 - \theta\varphi + \omega^{i_0} = 0 \quad (\varphi_1 > 1 > \varphi_2 > 0)$
We remark that the last equation includes k=1 because $U_{i_0}$ satisfies the pattern.

**Step 4**
Our focus is now the interval $[0, i_0]$.
$$U_1 = qzU_2$$
$$U_{i_0-1} = pzU_{i_0-2} + qz(1-s)U_{i_0}$$

so:
$$b_1 = -a_1 = \frac{(1-s)U_{i_0}}{\tau_2^{i_0} - \tau_1^{i_0}}$$

Using $U_0 = qzU_1$ we get: $\quad (\tau_2^{i_0} - \tau_1^{i_0}) U_0 = q(1-s)z(\tau_2 - \tau_1)U_{i_0}$

**Step 5**
The last step uses the starting point:
$$U_{i_0} = 1 + pzU_{i_0-1} + qzU_{i_0+1}$$
and after some calculation, using the preceding steps, we get:

$$U_{ki_0} = \frac{(\tau_2^{i_0} - \tau_1^{i_0})\varphi_2^k}{q(1-s)z(\tau_2-\tau_1)\omega^{i_0}} \quad (k \geq 1)$$

$$U_0 = \omega^{-i_0}\varphi_2$$

**CASE [z=1 and $\omega \neq 1$]**
We can follow the 5 steps above, the only difference is : $\tau_1 = \max(1, \omega), \tau_2 = \min(1, \omega)$.

CASE [z=1 and $\omega = 1$]

The strategy is the same; we give the result of each step.
Step 1
Because of $\tau_1 = \tau_2 = 1$ we have:
$$U_{ki_0+n} = a_{k+1}n + b_{k+1} \quad (k \geq 0, \quad n = 1,2,\dots,i_0-1)$$
Step 2
$$a_{k+1} = \frac{(1-s)[U_{(k+1)i_0}-U_{ki_0}]}{i_0} \quad (k>0)$$
$$b_{k+1} = (1-s)U_{ki_0} \quad (k>0)$$

Step 3
$$\varphi^2 - \theta\varphi + 1 = 0 \quad (\varphi_1 > 1 > \varphi_2 > 0)$$
$$\theta = 2\left(\frac{i_0}{1-s} + 1 - i_0\right)$$
$$U_{ki_0} = C_2\varphi_2^k \quad (k \geq 1),$$

Step 4
$$a_1 = \frac{(1-s)U_{i_0}}{i_0}; \quad b_1 = 0$$
Step 5
$$U_0 = U_0(1,1,s) = \varphi_2$$
$$U_{ki_0} = U_{ki_0}(1,1,s) = \frac{2i_0\varphi_2^k}{1-s} \quad (k \geq 1)$$

Remark: The same result can be obtained by applying l'Hospital's rule in the case 0<z<1

**Theorem 2**
The moment generating functions V for player B on the set $\{ki_0|\ k \geq 0\}$ are:
$$V_{i_0} = \frac{U_{i_0}-1}{1-s}$$
$$V_{ki_0} = \frac{U_{ki_0}}{1-s} \quad (k \neq 1)$$
Proof

We use the following notation:

$p_{i_0,ki_0}^{[X](m)}$=P(player X arrives in $ki_0$ in m steps before absorption has taken place, when starting in $i_0$).

The difference of player A and B is the first step: player A starts at t=0 with probabilities p(1-s) and q(1-s), while B is starting with p an q. The difference is a factor (1-s).

$$U_{ki_0} = \sum_{m=0}^{\infty} p_{i_0,ki_0}^{[A](m)} z^m = \delta(k,1) + \sum_{m=1}^{\infty} p_{i_0,ki_0}^{[A](m)} z^m =$$
$$\delta(k,1)+(1-s)\sum_{m=1}^{\infty} p_{i_0,ki_0}^{[B](m)} z^m = \delta(k,1) + (1-s)\,V_{ki_0}$$

**Theorem 3**
The moment generating functions W for player C on the set $\{ki_0|\ k \geq 0\}$ are:

CASE [0<z<1] $\qquad W_0 = W_0(z, \omega, s) = \frac{1}{\tau_1^{i_0} + \tau_2^{i_0} - \varphi_2}$

$$W_{i_0}(z, \omega, s) = \frac{(\tau_2^{i_0} - \tau_1^{i_0})}{qz(\tau_2 - \tau_1)(\tau_1^{i_0} + \tau_2^{i_0} - \varphi_2)}$$

$$W_{ki_0}(z, \omega, s) = \frac{(\tau_2^{i_0} - \tau_1^{i_0})\varphi_2^{k-1}}{q(1-s)z(\tau_2 - \tau_1)(\tau_1^{i_0} + \tau_2^{i_0} - \varphi_2)} \quad (k \geq 2)$$

where $\varphi_1, \varphi_2$ are solutions of $\varphi^2 - \theta\varphi + \omega^{i_0} = 0 \quad (\varphi_1 > 1 > \varphi_2 > 0)$

$$\theta = \theta(z, \omega, s) = \frac{\frac{\tau_2^{i_0} - \tau_1^{i_0}}{1-s} - 2pz(\tau_2^{i_0-1} - \tau_1^{i_0-1})}{qz(\tau_2 - \tau_1)}$$

and $\tau_1, \tau_2$ are solutions of $\qquad qz\tau^2 - \tau + pz = 0 \ (\tau_1 > 1 > \tau_2 > 0)$

CASE [z=1 and $\omega \neq 1$] $\qquad \tau_1 = \max(1, \omega)\ , \tau_2 = \min(1, \omega).$

CASE [z=1 and $\omega = 1$] $\qquad W_0 = \frac{1}{2 - \varphi_2}$

$$W_{i_0} = \frac{2i_0}{2 - \varphi_2}$$

$$W_{ki_0} = \frac{2i_0 \varphi_2^{k-1}}{(1-s)(2 - \varphi_2)} \quad (k \geq 2)$$

$$\varphi^2 - \theta\varphi + 1 = 0 \quad (\varphi_1 > 1 > \varphi_2 > 0)$$

$$\theta = \theta(1, 1, s) = 2\left(\frac{i_0}{1-s} + 1 - i_0\right)$$

Proof

First we prove the results of CASE [0<z<1] in 8 steps:

Step 1
As in proof of theorem1: $\qquad W_{ki_0+n} = a_{k+1}\tau_1^n + b_{k+1}\tau_2^n \qquad (k \geq 0, \quad n = 1,2, \ldots, i_0 - 1)$

Step 2
As in theorem 1, now with $k \geq 2$:

$$a_{k+1} = \frac{(1-s)\left[\tau_2^{i_0} W_{ki_0} - W_{(k+1)i_0}\right]}{\tau_2^{i_0} - \tau_1^{i_0}}$$

$$b_{k+1} = \frac{(1-s)\left[W_{(k+1)i_0} - \tau_1^{i_0} W_{ki_0}\right]}{\tau_2^{i_0} - \tau_1^{i_0}}$$

Step 3
The same result as in theorem 1, but now $k \geq 2$: $W_{ki_0} = C_2 \varphi_2^k$

**Step 4**
Our focus is the interval $[0, i_0]$.

$$W_1 = qzW_2; \quad W_{i_0-1} = pzW_{i_0-2} + qzW_{i_0}$$

$$b_1 = -a_1 = \frac{W_{i_0}}{\tau_2^{i_0} - \tau_1^{i_0}}$$

$$W_0 = qzW_1$$

$$(\tau_2^{i_0} - \tau_1^{i_0})W_0 = qz(\tau_2 - \tau_1)W_{i_0}$$

**Step 5**
Focus on $[i_0, 2i_0]$.

$$W_{i_0+1} = pzW_{i_0} + qzW_{i_0+2}$$
$$W_{2i_0-1} = pzW_{2i_0-2} + qz(1-s)W_{2i_0}$$

$$a_2 = \frac{\left[\tau_2^{i_0}W_{i_0} - (1-s)W_{2i_0}\right]}{\tau_2^{i_0} - \tau_1^{i_0}}$$

$$b_2 = \frac{\left[(1-s)W_{2i_0} - \tau_1^{i_0}W_{i_0}\right]}{\tau_2^{i_0} - \tau_1^{i_0}}$$

**Step 6**
This step uses the starting point:

$$W_{i_0} = 1 + pzW_{i_0-1} + qzW_{i_0+1}$$

$$qz(\tau_2 - \tau_1)\left[(\tau_2^{i_0} + \tau_1^{i_0})W_{i_0} - (1-s)W_{2i_0}\right] = (\tau_2^{i_0} - \tau_1^{i_0})$$

**Step 7**
Connection in $2i_0$:

$$W_{2i_0} = pzW_{2i_0-1} + qzW_{2i_0+1}$$

Using $\varphi^2 - \theta\varphi + \omega^{i_0} = 0$ we get:

$$\omega^{i_0}W_{i_0} - (1-s)\varphi_1 W_{2i_0} = 0$$

**Step 8**
Combining previous steps:

$$W_0 = \frac{1}{\tau_1^{i_0} + \tau_2^{i_0} - \varphi_2}$$

$$W_{i_0} = \frac{(\tau_2^{i_0} - \tau_1^{i_0})}{qz(\tau_2 - \tau_1)(\tau_1^{i_0} + \tau_2^{i_0} - \varphi_2)}$$

$$W_{ki_0} = \frac{(\tau_2^{i_0} - \tau_1^{i_0})\varphi_2^{k-1}}{q(1-s)z(\tau_2 - \tau_1)(\tau_1^{i_0} + \tau_2^{i_0} - \varphi_2)} \quad (k \geq 2)$$

**CASE [z=1 and $\omega \neq 1$]**
We can follow the 8 steps above, the only difference is : $\tau_1 = \max(1, \omega), \tau_2 = \min(1, \omega)$.

CASE [z=1 and $\omega = 1$]
The strategy is the same; we give the result of each step.

Step 1
Because of $\tau_1 = \tau_2 = 1$ we have:
$$W_{ki_0+n} = a_{k+1}n + b_{k+1} \quad (k \geq 0, \quad n = 1,2,\ldots, i_0 - 1)$$

Step 2
$$a_{k+1} = \frac{(1-s)[W_{(k+1)i_0} - W_{ki_0}]}{i_0} \quad (k>1)$$
$$b_{k+1} = (1-s)W_{ki_0} \quad (k>1)$$

Step 3
$$\varphi^2 - \theta\varphi + 1 = 0 \quad (\varphi_1 > 1 > \varphi_2 > 0)$$
$$\theta = 2\left(\frac{i_0}{1-s} + 1 - i_0\right)$$
$$W_{ki_0} = C_2\varphi_2^k \quad (k \geq 1),$$

Step 4
$$a_1 = \frac{W_{i_0}}{i_0} \; ; \quad b_1 = 0$$

Step 5
$$a_2 = \frac{[(1-s)W_{2i_0} - W_{i_0}]}{i_0} \; ; \quad b_2 = W_{i_0}$$

Step 6
$$2W_{i_0} - (1-s)W_{2i_0} = 2i_0$$

Step 7
$$W_{i_0} - (1-s)\varphi_1 W_{2i_0} = 0$$

Step 8
$$W_0 = \frac{1}{2-\varphi_2}$$
$$W_{i_0} = \frac{2i_0}{2-\varphi_2}$$
$$W_{ki_0} = \frac{2i_0\varphi_2^{k-1}}{(1-s)(2-\varphi_2)} \quad (k \geq 2)$$

Remark: Again the last results can be obtained by applying l'Hospital's rule in the CASE [0<z<1].

## 3. Absorption probabilities

It is easy to obtain the absorption probabilities in the absorption barrier 0 and in the mfb's.

We use the following notation: $P_X(Y) = $P(absorption in state Y for player X).

**Theorem 4**

CASE $\omega \neq 1$

Absorption probabilities player B:

$$P_B(0) = V_0(1, \omega, s) = \frac{\omega^{-i_0}\varphi_2}{1-s}$$

$$P_B(i_0) = sV_{i_0}(1, \omega, s) = s\left[\frac{U_{i_0}(1,\omega,s)-1}{1-s}\right] = \left[\frac{s}{1-s}\right]\left[\frac{(\omega^{-i_0}-1)\varphi_2}{(1-s)(q-p)} - 1\right]$$

$$P_B(ki_0) = sV_{ki_0}(1, \omega, s) = \frac{s(\omega^{-i_0}-1)\varphi_2^k}{(1-s)^2(q-p)} \quad (k \geq 2)$$

Player C:

$$P_C(0) = W_0(1, \omega, s) = \frac{1}{1+\omega^{i_0}-\varphi_2}$$

$$P_C(ki_0) = sW_{ki_0}(1, \omega, s) = \frac{s(1-\omega^{i_0})\varphi_2^{k-1}}{(1-s)(q-p)(1+\omega^{i_0}-\varphi_2)} \quad (k \geq 2)$$

CASE $\omega = 1$

Player B:

$$P_B(0) = V_0(1,1,s) = \frac{\varphi_2}{1-s}$$

$$P_B(i_0) = sV_{i_0}(1,1,s) = s\left[\frac{U_{i_0}(1,\omega,s)-1}{1-s}\right] = \left[\frac{s}{1-s}\right]\left[\frac{2i_0\varphi_2}{(1-s)} - 1\right]$$

$$P_B(ki_0) = sV_{ki_0}(1,1,s) = \frac{2i_0 s\varphi_2^k}{(1-s)^2} \quad (k \geq 2)$$

Player C:

$$P_C(0) = W_0(1,1,s) = \frac{1}{2-\varphi_2}$$

$$P_C(ki_0) = sW_{ki_0}(1,1,s) = \frac{2i_0 s\varphi_2^{k-1}}{(1-s)(2-\varphi_2)} \quad (k \geq 2)$$

We notice that there is a constant ratio between the B and C absorption probabilities:

CASE $\omega \neq 1$

$$\frac{P_B(0)}{P_C(0)} = \frac{P_B(ki_0)}{P_C(ki_0)} = \frac{\varphi_2(1 + \omega^{i_0} - \varphi_2)}{(1-s)\omega^{i_0}} < 1 \quad (k \geq 2)$$

CASE $\omega = 1$

$$\frac{P_B(0)}{P_C(0)} = \frac{P_B(ki_0)}{P_C(ki_0)} = \frac{\varphi_2(2 - \varphi_2)}{(1-s)} < 1 \quad (k \geq 2)$$

## 4. Mean absorption time

4.1 Mean time before absorption in any mfb

We define the mean time until absorption in any mfb for player X, when starting in i as:

$$m_i^{[X]} = \sum_{n=0}^{\infty} \sum_{k=0}^{\infty} n p_{i,ki_0}^{[X](n)} s_k$$

where $s_k$ is the (immediate) absorption probability in state $ki_0$ $(k \geq 0)$.
We have:

$$m_i^{[X]} = \sum_{n=1}^{\infty} \sum_{k=0}^{\infty} n p_{i,ki_0}^{[X](n)} s_k = \sum_{n=1}^{\infty} \sum_{k=0}^{\infty} (n-1) p_{i,ki_0}^{[X](n)} s_k + \sum_{n=1}^{\infty} \sum_{k=0}^{\infty} p_{i,ki_0}^{[X](n)} s_k =$$

$$\sum_{n=1}^{\infty} \sum_{k=0}^{\infty} (n-1) \left[ p p_{i+1,ki_0}^{[X](n-1)} + q p_{i-1,ki_0}^{[X](n-1)} \right] s_k + 1 = p m_{i+1}^{[X]} + q m_{i-1}^{[X]} + 1$$

**Theorem 5**
The mean absorption time in any mfb for player A is:
CASE $\omega \neq 1$
$$m_{i_0}^{[A]} = i_0 \frac{(1-s)}{s} \left(1 - \varphi_2 \omega^{-i_0}\right)$$
CASE $\omega = 1$
$$m_{i_0}^{[A]} = i_0 \frac{(1-s)}{s} (1 - \varphi_2)$$

Proof
CASE $\omega \neq 1$

Step 1
We describe the behaviour of $m_i^{[A]}$ between the multiple function barriers:
$$m_{ki_0+n}^{[A]} = p m_{ki_0+n+1}^{[A]} + q m_{ki_0+n-1}^{[A]} + 1$$

where at first glance $k \geq 0$, $n = 2,3,\ldots,i_0 - 2$, but n=1 and n=$i_0 - 1$ also satisfy the pattern.

We get:
$$m_{ki_0+n}^{[A]} = \frac{n}{q-p} + \omega^{-n} a_{k+1} + b_{k+1} \quad (k \geq 0, \ n = 1,2,\ldots,i_0 - 1)$$

Step 2
We express the constants of step 1 in the m values of the mfb's.
On the interval $[ki_0, (k+1)i_0]$, $k \geq 0$ we have

$$m^{[A]}_{ki_0+1} = pm^{[A]}_{ki_0+2} + qm^{[A]}_{ki_0} + 1$$

$$m^{[A]}_{(k+1)i_0-1} = pm^{[A]}_{(k+1)i_0} + qm^{[A]}_{(k+1)i_0-2} + 1$$

$$a_{k+1} = \frac{m^{[A]}_{ki_0} - m^{[A]}_{(k+1)i_0} - \frac{i_0}{p-q}}{1-\omega^{-i_0}}$$

$$b_{k+1} = \frac{m^{[A]}_{(k+1)i_0} + \frac{i_0}{p-q} - \omega^{-i_0} m^{[A]}_{ki_0}}{1-\omega^{-i_0}}$$

where $m^{[A]}_0 = 0$.

Step 3:
Focussing on the mfb in $ki_0$ ($k \geq 1$) we get:
$$m^{[A]}_{ki_0} = p(1-s)m^{[A]}_{ki_0+1} + q(1-s)m^{[A]}_{ki_0-1} + (1-s)$$

and after some calculation:

$$\omega^{i_0} m^{[A]}_{(k+1)i_0} - \theta m^{[A]}_{ki_0} + m^{[A]}_{(k-1)i_0} = \frac{i_0(1-\omega^{i_0})}{p-q} \quad (k \geq 1)$$

so we have:
$$m^{[A]}_{ki_0} = i_0 \frac{(1-s)}{s} + C_1 \varphi_1^{-k} \quad (k \geq 0),$$

and, using $m^{[A]}_0 = 0$:

$$m^{[A]}_{i_0} = i_0 \frac{(1-s)}{s}(1 - \varphi_1^{-1})$$

CASE $\omega = 1$

The strategy is the same; now starting with $m^{[A]}_{ki_0+n} = a_{k+1} + nb_{k+1} - n^2$.

**Theorem 6**
The mean absorption time in any mfb for player B is:
$$m^{[B]}_{i_0} = i_0(1 - \varphi_2 \omega^{-i_0})$$
Proof
The only difference of player A and B is the first step: player A starts at t=0 with probabilities p(1-s) and q(1-s), while B is starting with p an q. The difference in $m_{i_0}$ is a factor (1-s):

$$m^{[A]}_{i_0} = \sum_{n=1}^{\infty} \sum_{k=0}^{\infty} np^{[A](n)}_{i_0,ki_0} s_k = \sum_{n=1}^{\infty} \sum_{k=0}^{\infty} n(1-s) p^{[B](n)}_{i_0,ki_0} s_k = (1-s)m^{[B]}_{i_0}$$

**Theorem 7**
The mean absorption time in any mfb for player C is:
CASE $\omega \neq 1$

$$m_{i_0}^{[C]} = \frac{i_0\left\{\left(\frac{1-s}{s}\right)(1-\varphi_1^{-1}) + \left[\frac{1-\omega^{-i_0}}{p-q}\right]\right\}}{1+\omega^{-i_0}-\varphi_1^{-1}}$$

CASE $\omega = 1$

$$m_{i_0}^{[C]} = \frac{i_0\left[2i_0 + \left(\frac{1-s}{s}\right)(1-\varphi_2)\right]}{2-\varphi_2}$$

Proof

CASE $\omega \neq 1$

Steps 1,2 as in theorem 5

Step 3, now with $k \geq 1$: $\quad m_{ki_0}^{[C]} = i_0\left(\frac{1-s}{s}\right) + C_1\varphi_1^{-k} \quad (k \geq 1)$

Step 4 $\quad m_{i_0}^{[C]} = pm_{i_0+1}^{[C]} + qm_{i_0-1}^{[C]} + 1$

$$m_{2i_0}^{[C]} - (1+\omega^{-i_0})m_{i_0}^{[C]} = \frac{i_0(1-\omega^{-i_0})}{q-p}$$

and using the result of step 3 we get:

$$m_{i_0}^{[C]} = \frac{i_0\left\{\left(\frac{1-s}{s}\right)(1-\varphi_1^{-1}) + \left[\frac{1-\omega^{-i_0}}{p-q}\right]\right\}}{1+\omega^{-i_0}-\varphi_1^{-1}}$$

CASE $\omega = 1$

We give the results:
Step 1 $\quad m_{ki_0+n}^{[C]} = a_{k+1} + nb_{k+1} - n^2 \quad (k \geq 0)$

Step 2 $\quad a_{k+1} = m_{ki_0}^{[C]} \quad (k \geq 0)$

$$b_{k+1} = i_0 + \frac{m_{(k+1)i_0}^{[C]} - m_{ki_0}^{[C]}}{i_0} \quad (k \geq 0)$$

Step 3
$$m_{(k+1)i_0}^{[C]} - \theta(1,1,s)m_{ki_0}^{[C]} + m_{(k-1)i_0}^{[A]} = -2i_0^2 \quad (k \geq 1)$$

$$m_{ki_0}^{[A]} = i_0\left(\frac{1-s}{s}\right) + C_2\varphi_1^{-k} \quad (k \geq 1)$$

Step 4 $\quad m_{2i_0}^{[C]} - 2m_{i_0}^{[C]} = -2i_0^2$

$$m_{i_0}^{[C]} = \frac{i_0\left[2i_0 + \left(\frac{1-s}{s}\right)(1-\varphi_2)\right]}{2-\varphi_2}$$

## 4.2 Mean time before absorption in a specific mfb

We define:
$E\left[T_{ki_0}^{[X]}\right]$ = mean time until absorption in mfb $ki_0$, player X, when starting in $i_0$.

We have, if $\omega \neq 1$:

$$E\left[T_{ki_0}^{[X]}\right] = \sum_{n=0}^{\infty} n p_{i_0,ki_0}^{[X](n)} \cdot s_k = s_k \cdot \left(\frac{dM_X}{dz}\right)_{z=1}$$

where $s_k$ is the (immediate) absorption probability in state $ki_0$ ($k \geq 0$) and $M_X$ is the moment generating function of player X.

The connection with section 4.1 is:

$$m_{i_0}^{[X]} = \sum_{n=0}^{\infty}\sum_{k=0}^{\infty} n p_{i_0,ki_0}^{[X](n)} s_k = \sum_{k=0}^{\infty}\sum_{n=0}^{\infty} n p_{i_0,ki_0}^{[X](n)} s_k = \sum_{k=0}^{\infty} E\left[T_{ki_0}^{[X]}\right]$$

We discuss the CASE $\omega \neq 1$.
From theorems 1,2 and 3 we see that the mean time before absorption in the case $\omega \neq 1$ is closely related to differentiating the function $\varphi_2$.

Implicit differentiation of $qz\tau^2 - \tau + pz = 0$ gives

$$\frac{d\tau_i}{dz} = (-1)^i h(z) z^{-1} \tau_i \quad (i=1,2)$$

with
$$h(z) = (1 - 4pqz^2)^{-\frac{1}{2}}$$

$$\frac{dh(z)}{dz} = 4pqzh^3(z)$$

$$h(1) = \frac{1}{|p-q|}.$$

Implicit differentiation of $\varphi^2 - \theta\varphi + \omega^{i_0} = 0$ gives

$$\frac{d\varphi_i}{dz} = \left[\frac{\varphi_i}{2\varphi_i - \theta}\right]\frac{d\theta}{dz} = (-1)^i \left[\frac{\varphi_i}{\varphi_2 - \varphi_1}\right]\frac{d\theta}{dz} \quad (i=1,2)$$

with
$$\theta(z,\omega,s) = \frac{\frac{\tau_2^{i_0} - \tau_1^{i_0}}{1-s} - 2pz(\tau_2^{i_0-1} - \tau_1^{i_0-1})}{qz(\tau_2 - \tau_1)}$$

and
$$\left(\frac{d\theta}{dz}\right)_{z=1} = \frac{4pq\theta(1,\omega,s) - \frac{i_0(1+\omega^{i_0})}{1-s} + 2p[(p-q)(1-\omega^{i_0-1}) + (i_0-1)(1+\omega^{i_0-1})]}{(p-q)^2}$$

**Theorem 8**

$$\left(\frac{d\varphi_i}{dz}\right)_{z=1} = (-1)^i \left[\frac{\varphi_i}{\varphi_2 - \varphi_1}\right]\left[\frac{4pq\theta(1,\omega,s) - \frac{i_0(1+\omega^{i_0})}{1-s} + 2p[(p-q)(1-\omega^{i_0-1}) + (i_0-1)(1+\omega^{i_0-1})]}{(p-q)^2}\right] \quad (i=1,2)$$

**Appendix**

In this appendix we study two special cases: s=0 and s=1.

*A1. The case s=0*

There are no mfb's. Player's A,B and C now have the same strategy: play until ruin.
If s=0 then $\varphi_1 = \max(1, \omega^{i_0})$, $\varphi_2 = \min(1, \omega^{i_0})$.
Using theorem 4, we get:

$$P_{A,B,C}(0) = \omega^{-i_0}\varphi_2 = \begin{cases} 1 & if\ \omega \leq 1 \\ \omega^{-i_0} & if\ \omega > 1 \end{cases}$$

Using theorem 8:

$$m_0^{[A,B,C]} = E\left[T_0^{[A,B,C]}\right] = \left(\frac{d\varphi_2}{dz}\right)_{z=1} = \begin{cases} \dfrac{i_0\omega^{i_0}}{q-p} & if\ \omega \leq 1 \\ \dfrac{i_0}{p-q} & if\ \omega > 1 \end{cases}$$

(well known results, see [1]).

*A2. The case s=1*

Player B
Now we can't use player A as a reference.
We first calculate the moment generating function V on $(0, i_0)$ when starting in $i_0 - 1$:

$$V_n = \frac{(\tau_2^n - \tau_1^n)}{qz(\tau_2^{i_0} - \tau_1^{i_0})}\ ;\ V_0 = qzV_1 = \frac{\tau_2 - \tau_1}{(\tau_2^{i_0} - \tau_1^{i_0})}\ ;\ V_{i_0} = pzV_{i_0-1} = \frac{\omega(\tau_2^{i_0-1} - \tau_1^{i_0-1})}{(\tau_2^{i_0} - \tau_1^{i_0})}$$

Next the calculate V on $(i_0, 2i_0)$ when starting in $i_0 + 1$:

$$V_{i_0+n} = \frac{(\tau_1^n \tau_2^{i_0} - \tau_1^{i_0}\tau_2^n)}{pz(\tau_2^{i_0} - \tau_1^{i_0})}\ (0 < n < i_0);\ V_{2i_0} = pzV_{2i_0-1} = \frac{\omega^{i_0-1}(\tau_2 - \tau_1)}{(\tau_2^{i_0} - \tau_1^{i_0})}\ ;\ V_{i_0} = qzV_{i_0+1} = \frac{(\tau_2^{i_0-1} - \tau_1^{i_0-1})}{(\tau_2^{i_0} - \tau_1^{i_0})}$$

Absorption in 0 can only be obtained by moving one step backward on t=0, so we have:

$$p_0^{[B]} = qV_0(1) = \frac{q-p}{1-\omega^{i_0}}$$
$$p_{i_0}^{[B]} = q(pV_{i_0-1}) + p(qV_{i_0+1}) = \frac{2p(1-\omega^{i_0-1})}{(1-\omega^{i_0})}$$
$$p_{2i_0}^{[B]} = pV_{2i_0-1}(1) = \frac{(q-p)\omega^{i_0}}{1-\omega^{i_0}}$$

If $\omega = 1$ then:

$$p_0^{[B]} = p_{2i_0}^{[B]} = \frac{1}{2i_0}$$

$$p_{i_0}^{[B]} = \frac{i_0 - 1}{i_0}$$

Define $m_i^{*[B]}$ as the mean absorption time in any mfb, starting in i, for t>0: state $i_0$ has become an absorption barrier.

We have:
$$m_{i_0}^{[B]} = pm_{i_0+1}^{*[B]} + qm_{i_0-1}^{*[B]} + 1 = 1 + p\left[\frac{1}{q-p} + \frac{i_0}{p(1-\omega^{-i_0})}\right] + q\left[\frac{1}{p-q} + \frac{i_0}{q(1-\omega^{i_0})}\right] = i_0$$

We define $T_{i,j}^{[B]}$ as the time needed for player B for absorption in state j, when starting in i.
We have:
$$E[T_{i_0,0}^{[B]}] = \sum np_{i_0,0}^{(n)} = \sum (n-1)p_{i_0,0}^{(n)} + \sum p_{i_0,0}^{(n)} = \sum (n-1)qp_{i_0-1,0}^{(n-1)} + \sum p_{i_0,0}^{(n)} = qE[T_{i_0-1,0}^{[B]}] + p_0^{[B]}$$

where
$$E[T_{i_0-1,0}^{[B]}] = \left(\frac{dV_0}{dz}\right)_{z=1}$$

so:
$$E[T_{i_0,0}^{[B]}] = \frac{1}{(q-p)(\omega^{i_0}-1)} + \frac{i_0(\omega^{i_0}+1)}{(\omega^{i_0}-1)^2} + p_0^{[B]}$$

$$E[T_{i_0,2i_0}^{[B]}] = \omega^{i_0}\left[\frac{1}{(q-p)(\omega^{i_0}-1)} + \frac{i_0(\omega^{i_0}+1)}{(\omega^{i_0}-1)^2}\right] + p_{2i_0}^{[B]}$$

$$E[T_{i_0,i_0}^{[B]}] = -\left[\frac{2p(\omega^{i_0-1}+1)}{(q-p)(\omega^{i_0}-1)} + \frac{4i_0\omega^{i_0}}{(\omega^{i_0}-1)^2}\right] + p_{i_0}^{[B]}$$

and
$$E[T_{i_0,0}^{[B]}] + E[T_{i_0,2i_0}^{[B]}] + E[T_{i_0,i_0}^{[B]}] = i_0 = m_{i_0}^{[B]}$$

Player C
Using the same methods as in our main paper:
$$p_0^{[C]} = \frac{1}{(1+\omega^{i_0})}$$

$$p_{2i_0}^{[C]} = \frac{\omega^{i_0}}{(1+\omega^{i_0})}$$

$$m_{i_0}^{[C]} = \frac{i_0(1-\omega^{i_0})}{(q-p)(1+\omega^{i_0})}$$

$$E\left[T_0^{[C]}\right] = \frac{i_0(1-\omega^{i_0})}{(q-p)(1+\omega^{i_0})^2}$$

$$E\left[T_{2i_0}^{[C]}\right] = \frac{i_0(1-\omega^{i_0})\omega^{i_0}}{(q-p)(1+\omega^{i_0})^2}$$

We can simply verify that
$$m_{i_0}^{[C]} = E\left[T_0^{[C]}\right] + E\left[T_{2i_0}^{[C]}\right]$$